\documentclass{amsart}
\usepackage{amsfonts,amssymb}

\newcommand{\w}{\omega}
\newcommand{\IN}{\mathbb N}
\newcommand{\e}{\varepsilon}
\newcommand{\Iso}{\mathrm{Iso}}
\newcommand{\Ra}{\Rightarrow}

\newcommand{\IR}{\mathbb R}
\newcommand{\IZ}{\mathbb Z}
\newcommand{\A}{\mathcal A}
\newcommand{\F}{\mathfrak{F}}
\newcommand{\St}{\mathcal{S}t}

\newtheorem{theorem}{Theorem}[section]
\newtheorem{corollary}[theorem]{Corollary}
\newtheorem{proposition}[theorem]{Proposition}
\newtheorem{problem}[theorem]{Problem}
\newtheorem{question}[theorem]{Question}
\newtheorem{example}[theorem]{Example}
\newtheorem{lemma}[theorem]{Lemma}
\theoremstyle{definition}
\newtheorem{definition}[theorem]{Definition}

\newtheorem{remark}[theorem]{Remark}

\title[Kaleidoscopical configurations]{Kaleidoscopical configurations\\ in $G$-spaces}
\author{T.~Banakh, O.~Petrenko, I.V.~Protasov, S.~Slobodianiuk}
\address{Faculty of Mechanics and Mathematics, Lviv University, Universytetska 1, Lviv 79000, Ukraine}
\email{t.o.banakh@gmail.com}
\address{Faculty of Cybernetics, Kyiv University, Volodymyrska 64, Kyiv 01033, Ukraine}
\email{opetrenko72@gmail.com; i.v.protasov@gmail.com; }
\address{Faculty of Mechanics and Mathematics, Kyiv University, Volodymyrska 64, Kyiv 01033, Ukraine}
\email{slobodianiuk@yandex.ru}

\begin{document}
\begin{abstract}
Let $G$ be a group and $X$ be a $G$-space with the action $G\times X\rightarrow X$, $(g,x)\mapsto gx$. A subset $F$ of $X$ is called a kaleidoscopical configuration if there exists a coloring $\chi:X\rightarrow C$ such that the restriction of $\chi$ on each subset $gF$, $g\in G$, is a bijection. We present a construction (called the splitting construction) of kaleidoscopical configurations in an arbitrary $G$-space, reduce the problem of characterization of kaleidoscopical configurations in a finite Abelian group $G$ to a factorization of $G$ into two subsets, and describe all kaleidoscopical configurations in isometrically homogeneous ultrametric spaces with finite distance scale. Also we construct $2^{\mathfrak c}$ many (unsplittable) kaleidoscopic configurations of cardinality $\mathfrak c$ in the Euclidean space $\IR^n$.
\end{abstract}
\maketitle

\section*{Introduction}
Let $X$ be a set and $\mathfrak{F}$ be a family of subsets of $X$
(the pair $(X,\mathfrak{F})$ is called a {\em hypergraph}).
Following \cite{BP}, we say that a coloring $\chi:X\to\kappa$ of $X$
(i.e. a surjective mapping of $X$ onto a cardinal $\kappa$) is
\begin{itemize}
\item $\mathfrak{F}$-{\em surjective} if the restriction $\chi|_F$ is surjective for all $F\in\mathfrak{F}$;
\item $\mathfrak{F}$-{\em injective} if $\chi|_F$ is injective for all $F\in\mathfrak{F}$;
\item $\mathfrak{F}$-{\em bijective} or $\mathfrak{F}$-{\em kaleidoscopical} if $\chi|_F$ is bijective for all $F\in\mathfrak{F}$.
\end{itemize}

A hypergraph $(X,\mathfrak{F})$ is called {\em kaleidoskopical} if
there exists an $\mathfrak{F}$-kaleidoscopical coloring
$\chi:X\to\kappa$. The adjactive "kaleidoscopical" appeared in
definition \cite{PP} of $s$-regular graph $\Gamma(V,E)$ (each vertex
$v\in V$ has degree $s$) admitting a vertex $(s+1)$-colloring such
that each unit ball $B(v,1)=\{u\in V:d(u,v)=1\}$ has the vertices of
all colors ($d$ is a path matric on $V$). These graphs can be
considered as a graph version of Hamming codes \cite{PO}.

We shall consider hypergraphs related to $G$-space. Let $G$ be a
group. A $G$-{\em space} is a set $X$ endowed with an action
$G\times X\to X$, $(g,x)\mapsto gx$. All $G$-spaces are suppose to
be transitive (for any $x,y\in X$ there exists $g\in G$ such that
$gx=y$). For a subset $A\subseteq X$, we put $G[A]=\{gA:g\in G\}$.

A subset $A\subseteq X$ is called a {\em kaleidoscopical
configuration} if the hypergraph $(X,G[A])$ is kaleidoscopical (in
words, if there exists a coloring $\chi:X\to |A|$ such that
$\chi|gA$ is bijective for every $g\in G$).

In Section~\ref{s:groups} we show that kaleidoscopical
configurations are tightly connected with classical combinatorial
theme {\em Transversality} and, in the case $X=G$ and (left) regular
action of $G$ on $G$, with factorization problem, well known in {\em
Factorization Theory of groups}, see \cite{Sz}, \cite{SS}.

In Section~\ref{s:split} we introduce and describe the
kaleidoscopical configurations (called splittable) which arise from
the chains of $G$-invariant equivalences (imprimitivities) on $X$.
If a $G$-space $X$ is primitive (the only $G$-invariant equivalences
on $X$ are $X\times X$ and $\Delta_X$) then the only splittable
configurations in $X$ are $X$ and the singletons.

In Section~\ref{s:prim} we prove that every kaleidoscopical
configuration in isometrically homogeneous metric space with finite
distance scale is splittable. For $n\geq 2$, we construct a plenty
of kaleidoscopical configurations of cardinality $\mathfrak{c}$ in
$\IR^n$. These configurations are non-splittable because $\IR^n$ is
isometrically primitive. We don't know whether there exists a finite
non-singleton or countable kaleidoscopical configurations in
$\IR^n$, $n\geq 2$.

In Section~\ref{s:Hajos} we  study the problem of splittability of
kaleidoscopic configurations in finite Abelian groups and
reformulate this problem in terms of the semi-Haj\'os property, see
\cite{Sz}, \cite{SS}.

\section{Transversality and factorization}\label{s:groups}

Let $(X,\mathfrak{F})$ be a hypergraph. A subset $T\subseteq X$ is
called an $\mathfrak{F}$-transversal if $|F\bigcap T|=1$ for each
$F\in\mathfrak{F}$.

\begin{proposition} A hypergraph $(X,\mathfrak{F})$ is
kaleidoscopical if and only if $X$ can be partitioned into
$\mathfrak{F}$-transversals.
\end{proposition}
\begin{proof}
For a kaleidoscopic hypergraph $(X,\mathfrak{F})$, let
$\chi:X\to\kappa$ be a kaleidoscopical coloring. Then
$\bigsqcup_{\alpha<\kappa}\chi^{-1}(\alpha)$ is a partition of $X$
into $\mathfrak{F}$-transversal.

On the other hand, if $\bigsqcup_{\alpha<\kappa}T_\alpha$ is a
partition of $X$ into $\mathfrak{F}$-transversal then the coloring
$\chi:X\to\kappa$ defined as $\chi(x)=\alpha\Leftrightarrow x\in
T_\alpha$ is kaleidoscopical.
\end{proof}

Let $X$ be a $G$-space $A$ be a kaleidoscopical configuration in
$X$. If $T$ is $G[A]$-transversal then $A$ is $G[T]$-transversal and
$gT$ is $G[A]$ transversal for each $g\in G$.

We say  that a kaleidoscopical configuration $A$ in $X$ is {\em
homogeneous} if there exist a $G[A]$-transversal $T$ and a subset
$H\subseteq X$ such that $X=\bigsqcup_{h\in H}hT$.

A subset $A$ of a group $G$ is defined to be {\em complemented} in
$G$ if there exists a subset $B\subseteq G$ such that the
multiplication mapping $\mu:A\times B\to G$, $(a,b)\mapsto ab$, is
bijective. Following \cite{SS}, we call the set $B$ a {\em
complementer factor} to $A$, and say that $G=AB$ is a factorization
of $G$. In this case, we have
$$G=\bigsqcup_{a\in A}aB=\bigsqcup_{b\in B}Ab.$$
A subset $A\subseteq G$ is called {\em doubly complemented} if there
are factorization $G=AB=BC$ for some subsets $B,C$ of $G$.
\begin{proposition}\label{factor} For two subsets $A,B$ of a group $G$ the following conditions are equivalent:
\begin{enumerate}
\item $B$ is $G[A]$-transversal;
\item $G=AB^{-1}$ is a factorization of $G$.
\end{enumerate}
\end{proposition}

\begin{proof} $(1)\Ra(2)$ For each $g\in G$, $g^{-1}A\cap B\ne\emptyset$, so $g\in AB^{-1}$. If $g=a_1b_1^{-1}=a_2b_2^{-1}$ for some $a_1,a_2\in A$, $b_1,b_2\in B$, then $g^{-1}a_1=b_1$ and $g^{-1}a_2^{-1}=b_2$ and by (1), $b_1=b_2$ and $a_1=a_2$, witnessing that $G=A B^{-1}$ is a factorization of $G$.
\smallskip

$(2)\Ra(1)$ Fix any $g\in G$. The inclusion $g^{-1}\in AB^{-1}$ implies $gA\cap B\ne\emptyset$. If $ga_1=b_1$ and $ga_2=b_2$ for some $a_1,a_2\in A$, $b_1,b_2\in B$, then $g^{-1}=a_1b_1^{-1}=a_2b_2^{-1}$ and by (2), $b_1=b_2$, witnessing that $|gA\cap B|=1$.
\end{proof}
\begin{corollary} Each kaleidoscopic configuration in group $G$
is complemented.
\end{corollary}
\begin{proof} Given a kaleidoscopical configuration $A\subset G$, fix a $A$-kaleidoscopical coloring $\chi:G\to C$. We choose a color $c\in C$, consider the monochrome class $B=\chi^{-1}(b)\subset G$ and observe that for every $g\in G$ \ $|gA\cap B|=1$ by the definition of $A$-kaleidoscopical coloring. By Proposition~\ref{factor}, $G=AB^{-1}$ is a factorization, so $A$ is complemented in $G$.
\end{proof}
\begin{proposition}\label{p1.4} A subset $A$ of a group $G$ is doubly
complemented if and only if $A$ is a homogeneous kaleidoscopical
configuration.
\end{proposition}
\begin{proof} Let $G=AB=BC$ be a factorization of $G$. By
proposition~\ref{factor}, $B^{-1}$ is a $G[A]$-transversal. Since
$G=\bigsqcup_{c\in C}c^{-1}B$, we conclude that $A$ is a homogeneous
kaleidoscopical configuration.

Let $A$ be a homogeneous kaleidoscopical configuration. We choose a
$G[A]$-transversal $T$ and a subset $H\subseteq G$ such that
$G=\bigsqcup_{h\in H}hT$. By proposition~\ref{factor}, $G=AT^{-1}$.
Since $G=\bigsqcup_{h\in H}hT$, $G=T^{-1}H^{-1}$ is a factorization.
Hence, $A$ is doubly complemented.
\end{proof}
\begin{corollary} For a subset $A$ of an Abelian group $G$, the
following statements are equivalent:
\begin{enumerate}
\item $A$ is complemented;
\item $A$ is a kaleidoscopical configuration;
\item  $A$ is a homogeneous kaleidoscopical configuration.
\end{enumerate}
\end{corollary}
\begin{question} Is each complemented subset of a (finite) group
kaleidoscopical?
\end{question}
\begin{proposition}\label{section} Let $X$ be a $G$-space, $x\in X$, $G_x=\{g\in
G:gx=x\}$, $\gamma_x: G\to X$, $\gamma_x(g)=gx$, $s:X\to G$ be a
section of $\gamma_x$. Let $A$ be a subset of $X$, $T$ be
$G[A]$-transversal. Then
\begin{enumerate}
\item $s(T)$ is a $G[\gamma_x^{-1}(A)]$-transversal;
\item $|G|=|G_x||A||T|$.
\end{enumerate}
\end{proposition}
\begin{proof} The statement $(1)$ is evident. The statement $(2)$
follows from $(1)$ and proposition~\ref{factor}.
\end{proof}
\begin{corollary}\label{c1.8} Let $A$ be a kaleidoscopical configuration in a
finite $G$-space $X$ with a kaleidoscopical coloring  $\chi:G\to k$.
Then $\chi^{-1}(0)=\dots=\chi^{-1}(k-1)$ and
$|X|=|A||\chi^{-1}(0)|$.
\end{corollary}
\begin{proof} We may suppose that $G$ is a subgroup of the group of
all permutations of $X$ so $G$ is finite. Since $|G|=|X||G_x|$, we
can apply proposition~\ref{section}(2).
\end{proof}

\begin{proposition} Let $\kappa$ be an infinite cardinal, $(X,\mathfrak{F})$
be a hypergraph such that $|\mathfrak{F}|=\kappa$ and $|F|=\kappa$
for each $F\in\mathfrak{F}$. If $|F\cap F'|<cf\kappa$ for all
distinct $F,F'\in\mathfrak{F}$ then there are a disjoint family
$\mathfrak{T}$ of $\mathfrak{F}$-transversals such that
$|\mathfrak{T}|=\kappa$, $|T|=\kappa$ for each $T\in\mathfrak{T}$
\end{proposition}
\begin{proof} We enumerate $\mathfrak{F}=\{F_\alpha:\alpha<\kappa\}$
and choose inductively the subsets $\{V_\alpha\subset
F_\alpha:\alpha<\kappa\}$ such that the family $\{F_\alpha\setminus
V_\alpha:\alpha<\kappa\}$ is disjoint and $|F_\alpha\setminus
V_\alpha|=\kappa$ for each $\alpha<\kappa$. Let $F_\alpha\setminus
V_\alpha=\{t_{\alpha\beta}:\beta<\kappa\}$,
$T_\beta=\{t_{\alpha\beta}:\alpha<\kappa\}$. Then
$\mathfrak{T}=\{T_\beta:\beta<\kappa\}$ is the desired family.
\end{proof}
For a hypergraph $(X,\mathfrak{F})$, $x\in X$ and $A\subseteq X$, we
put
$$St(x,\mathfrak(F))=\bigcup\{F\in\mathfrak{F}:x\in F\},$$
$$St(A,\mathfrak{F})=\bigcup\{St(a,F):a\in A\}.$$
\begin{proposition}\label{p1.10} A hypergraph $(X,\mathfrak{F})$ is
kaleidoscopical provided that, for some infinite cardinal $\kappa$,
the following two conditions are satisfied:
\begin{enumerate}
\item $\mathfrak{F}\leq\kappa$ and $|F|=\kappa$ for each $F\in\mathfrak{F}$;
\item for any subfamily $\mathfrak{A}\subset\mathfrak{F}$ of cardinality $|\mathfrak{A}|<\kappa$ and any subset $B\subset X\setminus(\bigcup\mathfrak{A})$ of cardinality $|B|<\kappa$ the intersection $St(B,\mathfrak{F})\cap(\bigcup\mathfrak{A})$ has cardinality less then $\kappa$.
\end{enumerate}
\end{proposition}
\begin{proof} Let $\lambda=|\F|$ and $\F=\{F_\alpha:\alpha<\lambda\}$ be an injective enumeration of $\F$. By induction we shall construct a transfinite sequence $(\chi_\alpha:F_\alpha\to\kappa)_{\alpha<\lambda}$ of bijective colorings such that for any ordinals $\alpha<\beta<\lambda$
\begin{enumerate}
\item[$(1_{\alpha\beta})$] the colorings $\chi_\alpha$ and $\chi_\beta$ coincide on $F_\alpha\cap F_\beta$;
\item[$(2_{\alpha\beta})$] no distinct points $a\in F_\alpha$ and $b\in F_\beta$ with $\chi_\alpha(a)=\chi_\beta(b)$ lie in some hyperedge $F\in\F$.
\end{enumerate}

Assume that for some ordinal $\gamma<\lambda$ we have constructed a
sequence of colorings $(\chi_{\alpha})_{\alpha<\gamma}$ satisfying
the conditions $(1_{\alpha\beta})$ and $(2_{\alpha\beta})$ for all
$\alpha<\beta<\gamma$.

Let us define a bijective coloring $\chi_\gamma:F_\gamma\to\kappa$.
First we show that the union
$$F_\gamma'=\bigcup_{\alpha<\gamma}F_\gamma\cap F_\alpha$$ has
cardinality $|F'_\gamma|<\kappa$. Observe that for each
$\alpha<\gamma$ we get $F_\alpha\not\subset F_\gamma$. Assuming
conversely that $F_\alpha\subsetneq F_\gamma$ and taking any point
$v\in F_\gamma\setminus F_\alpha$ we conclude that the intersection
$F_\alpha\cap\St(v,\F)\supset F_\alpha\cap F_\gamma=F_\alpha$ has
cardinality $\ge\kappa$, which contradicts the condition (2) of the
theorem.

Therefore, for each $\alpha<\gamma$ we can choose a point
$v_\alpha\in F_\alpha\setminus F_\gamma$. Then for the set
$B=\{v_\alpha:\alpha<\gamma\}$ the set $F_\gamma'\subset
F_\gamma\cap \St(B,\F)$ has cardinality
$|F_\gamma'|\le|F_\gamma\cap\St(A,\F)|<\kappa$ according to (2).

For every point $x\in F_\gamma\setminus F_\gamma'$ and every ordinal
$\alpha<\gamma$ consider the sets $\St(x,\F)\cap F_\alpha$ and
$C_\alpha(x)=\chi_\alpha(\St(x,\F)\cap F_\alpha)\subset\kappa$. The
condition (2) implies that the set
$C(x)=\bigcup_{\alpha<\gamma}C_\alpha(x)$ has cardinality
$|C(x)|<\kappa$.

Let $\prec$ be any well-order on the set $F_\gamma$ such that
$F_\gamma'$ coincides the initial segment $\{x\in F_\gamma:x<y\}$
for some point $y\in F_\gamma$. Consider the coloring
$\chi_\gamma:F_\gamma\to\kappa$ defined by
$\chi_\gamma(x)=\chi_\alpha(x)$ if $x\in F_\gamma\cap F_\alpha$ for
some $\alpha<\gamma$ and
$$\chi_\gamma(x)=\min\,\kappa\setminus(C(x)\cup\{\chi(y):y\prec
x\})$$ if $x\in F_\gamma\setminus F_\gamma'$.

Let us show that the coloring $\chi_\gamma:F_\gamma\to \kappa$ is
bijective. The injectivity of $\chi_\gamma$ follows from the
definition of $\chi_\gamma$ and the conditions $(2_{\alpha\beta})$,
$\alpha<\beta<\gamma$.

The surjectivity of $\chi_\gamma$ will follow as soon as we check
that for each color $c\in\kappa\setminus\chi_\gamma(F_\gamma')$ the
set $F_\gamma(c)=\{x\in F_\gamma\setminus F_\gamma':c\in C(x)\}$ has
cardinality $<\kappa$. Observe that $c\in C(x)$ if and only if there
is $\alpha<\gamma$ and a point $a\in F_\alpha\setminus F_\gamma$
such that $\chi_\alpha(a)=c$ and $x\in \St(a,\F)$. The set
$A_c=\bigcup_{\alpha<\gamma}\chi_\alpha^{-1}(c)\setminus F_\gamma$
has size $|A_c|\le\gamma<\kappa$ and by the condition (2), the set
$F_\gamma(c)\subset F_\gamma\cap\St(A_c,\F)$ has cardinality
$<\kappa$. This completes the proof of the bijectivity of the
coloring $\chi_\gamma$.

The conditions $(1_{\alpha\gamma})$ and $(2_{\alpha,\gamma})$ for
all $\alpha<\gamma$ follow from the definition of the coloring
$\chi_\gamma$. This completes the inductive step of the construction
of the sequence $(\chi_\alpha)_{\alpha<\lambda}$.

After completing the inductive construction, let $\chi:V\to\kappa$
be any coloring such that $\chi|F_\alpha=\chi_\alpha$ for all
$\alpha<\lambda$. The conditions $(1_{\alpha\beta})$ guarantee that
the coloring $\chi$ is well-defined. The bijectivity of the
colorings $\chi_\alpha$, $\alpha<\lambda$, ensures the
kaleidoscopicity of the coloring $\chi$.
\end{proof}
We conclude this section with short discussion of possibilities of
transfering above notions and results of quasigroups.

We recall that a {\it quasigroup} is a set $X$ endowed with a binary
operation $*:X\times X\to X$ such that, for every $a,b\in X$, the
system of equations $a*x=b$, $y*a=b$ has a unique solution
$x=a\backslash b$, $y=b/a$ in $X$.

In an obvious way the notion of a kaleidoscopical configuration generalizes to quasigroup.

A subset $A$ of a quasigroup $X$ is called
\begin{itemize}
\item {\it kaleidoscopical} if there is a coloring $\chi:X\to C$ such that $\chi|_{x*A}:x*A\to C$ is bijective for all $x\in X$;
\item {\it complemented} if there is a subset $B\subset X$ such that the right division $\delta:B\times A\to X$, $\delta(b,a)=b/a$ is bijective;
\item {\it doubly complemented} if there exists a complemented subset  $B\subset X$ such that the multiplication $\mu:A\times B\to X$, $\mu(a,b)=a*b$,  is bijective;
\item {\it self-complemented} if the maps $\mu:A\times A\to X$, $\mu(x,y)=x*y$, and $\delta:A\times A\to X$, $\delta(x,y)=x/y$, are bijective.
\end{itemize}
It follows from the proof of proposition~\ref{factor} that each
kaleidoscopical subset in a semigroup is complemented. In contrast,
Proposition~\ref{p1.4} does not generalize to quasigroup.
\begin{example} There exists a quasigroup $X$ of order $|X|=9$ that contains a self-complemented subset $A\subset X$, which is not kaleidoscopical.
\end{example}

\begin{proof} It is well-known that finite quasigroups can be identified with
Latin squares, i.e., $n\times n$ matrices whose rows and columns are
permutations of the set $\{1,\dots,n\}$. For $r,s\le n$ an $(r\times
s)$-matrix $(x_{ij})$ is called a {\em partial Latin $(r\times
s)$-rectangle} if $x_{ij}\in\{1,2,\dots,n\}$ and $x_{lj}\ne
x_{ij}\ne x_{ik}$ for any $1\le i\ne l\le r$ and $1\le j\ne k\le s$.
By a result of Ryser \cite{Ryser}  (see also Lemma 1 in
\cite[p.214]{PB}) each partial latin $(r\times s)$-rectangle can be
completed to a Latin $(n\times n)$-square if and only if each number
$i\in\{1,\dots,n\}$ appears in the rectangle not less than $r+s-n$
times. This extension result allows us to find a quasigroups
operation on $X=\{1,\dots,9\}$ whose multiplication table has the
following first three columns:
\begin{table}[ht]
\begin{tabular}{c|ccc}
$*$&1&2&3\\
\hline
1&1&4&5\\
2&6&2&7\\
3&8&9&3\\
\hline
4&4&1&6\\
5&5&6&1\\
\hline
6&2&7&8\\
7&7&8&2\\
8&3&5&9\\
9&9&3&4\\
\hline
\end{tabular}
\end{table}

Looking at this table we can see that the set $A=\{1,2,3\}$ is self-complemented as $A*A=X=A/A$. Assuming that $A$ is kaleidoscopic, find a coloring $\chi:X\to A$ such that $\chi|x*A$ is bijective for each $x\in X$.
Since $1*A=\{1,4,5\}$ and $4*A=\{4,1,6\}$, the elements $5$ and $6$ have the same color, which is not possible as $5*A=\{5,6,1\}$ and $\chi|5*A$ is bijective.
\end{proof}

Corollary~\ref{c1.8} implies that the size $|K|$ of any
kaleidoscopic subset $K$ in a finite group $G$ divides the
cardinality $|G|$ of $G$. The same is true for any finite transitive
$G$-space $X$.

\section{Splitting}\label{s:split}

In this section we present a simple construction of kaleidoscopic configurations in arbitrary $G$-space, called the splitting construction. Kaleidoscopic subsets constructed in this way will be called splittable.

First we recall some definitions. A map $\varphi:X\to Y$ between $G$-spaces is called {\em equivariant} if $\varphi(gx)=g\,\varphi(x)$ for all $g\in G$ and  $x\in X$. It is easy to see that each equivariant map between transitive $G$-spaces is surjective and homogeneous.

A function $\varphi:X\to Y$ is defined to be {\em homogeneous} if it is $\kappa$-to-1 for some non-zero cardinal $\kappa$. The latter means that $|\varphi^{-1}(y)|=\kappa$ for all $y\in Y$.

\begin{proposition}\label{p2.1} Let $\kappa$ be a non-zero cardinal, $\pi:X\to Y$ be an $\kappa$-to-1 equivariant map between two $G$-spaces and $s:Y\to X$ be a section of $\varphi$. Let $K\subset Y$ be a kaleidoscopic subset and $\chi:Y\to C$ be an $K$-kaleidoscopic coloring. Then:
\begin{enumerate}
\item the preimage $\bar K=\pi^{-1}(K)$ is a kaleidoscopic configuration in $X$ with respect to any coloring $\bar \chi:X\to C\times \kappa$ such that for each $y\in Y$ the restriction $\bar\chi|\varphi^{-1}(y):\pi^{-1}(y)\to\{\chi(y)\}\times \kappa$ is bijective;
\item the image $\tilde K=s(K)$ is a kaleidoscopic configuration in $X$ with respect to the $\tilde K$-kaleidoscopic coloring $\tilde\chi=\chi\circ \pi:X\to C$.
\end{enumerate}
\end{proposition}

\begin{proof} 1. Given any element $g\in G$, we need to check that the restriction $\bar \chi|g\bar K:g\bar K\to C\times\kappa$ is bijective. To see that it is surjective, take any color $(c,\alpha)\in C\times\kappa$ and using the surjectivity of $\chi|gK:gK\to C$, find a point $y\in gK$ with $\chi(y)=c$. Since the restriction $\bar\chi|\pi^{-1}(y):\pi^{-1}(y)\to\{c\}\times\kappa$ is bijective, there is a point $x\in \pi^{-1}(y)\subset \pi^{-1}(gK)=g\bar K$ with $\bar\chi(x)=(c,\alpha)$, so $\bar\chi|g\bar K$ is surjective.

To see that it is injective, take any two distinct points $x,x'\in g\bar K$. If $\pi(x)=\pi(x')$, then for the point $y=\pi(x)=\pi(x')\in g\pi(K)=\pi(gK)$ the injectivity of the restriction $\bar \chi|\pi^{-1}(y)$ implies that $\bar\chi(x)\ne\bar\chi(x')$. If $\pi(x)\ne \pi(x')$, then the injectivity of $\chi|gK$ guarantees that $\chi(\pi(x))\ne\chi(\pi(x'))$ and then $\bar\chi(x)\ne\bar\chi(x')$ as $\bar\chi(x)\in\{\chi(\pi(x))\}\times\kappa$ and $\bar\chi(x)\in\{\chi(\pi(x'))\}\times\kappa$.
\smallskip

2. Given any element $g\in G$, we need to check that the restriction $\tilde \chi|g\tilde K:g\bar K\to C$ is bijective. To see that  $\tilde\chi|g\tilde K$ is surjective, observe that
$$\tilde \chi(g\tilde K)=\chi\circ \pi(g\tilde K)=\chi(g\,\pi(\tilde K))=\chi(gK)=C$$by the surjectivity of $\chi|gK:gK\to C$.

To see that $\tilde\chi|g\tilde K$ is injective, take any two distinct points $x,x'\in\tilde K=s(K)$ and observe $\pi(x)\ne \pi(x')$. Since $\pi$ is equivariant, $\pi(gx)=g\pi(x)\ne g\pi(x')=\pi(gx')$. Since $\pi(gx),\pi(gx')\in gK$ and $\chi|gK$ is injective, $\tilde \chi(x)=\chi(\pi(gx))\ne \chi(\pi(gx'))=\tilde\chi(\pi(gx'))$ are we are done.
\end{proof}

Iterating the constructions from Proposition~\ref{p2.1}, we get the so-called {\em splitting construction} of kaleidoscopical configurations.

\begin{proposition}\label{p2.2} Let $X_0\to X_{1}\to \cdots \to X_m$ be a sequence of $G$-spaces linked by homogeneous $G$-equivariant maps $\pi_i:X_i\to X_{i+1}$, $i<m$. Let $K_i\subset X_i$, $i\le m$, be subsets such that for every $i<m$ either the restriction $\pi_i|K_i:K_i\to K_{i+1}$ is bijective or else $K_i=\pi_i^{-1}(K_{i+1})$. If the set $K_m$ is kaleidoscopic in the $G$-space $X_m$, then for every $i\le m$ the set $K_i$ is kaleidoscopic in the $G$-space $X_i$.
\end{proposition}

\begin{proof} This proposition can be derived from Proposition~\ref{p2.1} by the reverse induction on $i\in\{m,m-1,\dots,0\}$.
\end{proof}

Proposition~\ref{p2.2} can be alternatively written in terms of invariant equivalence relations.

Given an equivalence relation $E\subset X\times X$ on a set $X$ let $X/E=\{[x]_E:x\in X\}$ be the quotient space consisting of the equivalence classes $[x]_E=\{y\in X:(x,y)\in E\}$, $x\in X$. Denote by $q_E:X\to X/E$, $q_E:x\mapsto [x]_E$,  the quotient map. For a subset $K\subset X$ let $K/E=\{[x]_E:x\in K\}\subset X/E$ and $[K]_E=\bigcup\limits_{x\in K}[x]_E\subset X$.

Let $E$ be an equivalence relation on a set $X$. A subset $K\subset X$ is  defined to be
\begin{itemize}
\item {\em $E$-parallel} \ if $K\cap [x]_E=[x]_E$ for all $x\in K$;
\item {\em $E$-orthogonal} \ if $K\cap [x]_E=\{x\}$ for all $x\in K$.
\end{itemize}
Given two equivalence relations $E\subset F$ on $X$ we can generalize these two notions defining $K\subset X$ to be
\begin{itemize}
\item {\em $F/E$-parallel} \ if $[K]_E\cap [x]_F=[x]_F$ for all $x\in K$;
\item {\em $F/E$-orthogonal} \ if $[K]_E\cap [x]_F=[x]_E$ for all $x\in K$.
\end{itemize}
Observe that a set $K\subset X$ is $E$-parallel ($E$-orthogonal) if and only if it is $E/\Delta_X$-parallel ($E/\Delta_X$-orthogonal). Here $\Delta_X=\{(x,x):x\in X\}$ stands for the smallest equivalence relation on $X$.


 An equivalence relation $E$ on a $G$-space $X$ is called
{\em $G$-invariant\/} if for each $(x,y)\in E$ and any $g\in G$ we get $(gx,gy)\in E$. For a $G$-invariant equivalence relation $E$ on $X$ the quotient space $X/E$ is a $G$-space under the induced action $$G\times X/E\to X/E,\;\;(g,[x]_E)\mapsto [gx]_E$$of the group $G$. In this case the quotient projection $q:X\to X/E$ is equivariant. $G$-Invariant equivalence relations on $G$-spaces are also called {\em imprimitivities}.

\begin{proposition}\label{p2.3} Let $\Delta_X=E_0\subset E_{1}\subset\dots\subset E_m$ be a sequence of $G$-invariant equivalence relations on a transitive $G$-space $X$. A subset $K\subset X$ is kaleidoscopic provided
\begin{enumerate}
\item the projection $K/E_m$ is kaleidoscopic in the $G$-space $X/E_m$;
\item for every $i<m$ the set $K$ is $E_{i+1}/E_i$-parallel or $E_{i+1}/E_i$-orthogonal.
\end{enumerate}
\end{proposition}

\begin{proof} For every $i\le m$ consider the $G$-space $X_i=X/E_i$ and the subset $K_i=K/E_i$ in $X_i$. Since $E_0=\Delta_X$, the space $X_0$ coincides with $X$. Next, for every $i<m$, consider the equivariant map $\pi_i:X_i\to X_{i+1}$, $\pi_i:[x]_{E_i}\mapsto [x]_{E_{i+1}}$. This map is homogeneous because of the transitivity of the $G$-space $X_i$.

We claim that the maps $\pi_i$ satisfy the requirements of
Proposition~\ref{p2.2}. Indeed, if $K$ is $E_{i+1}/E_i$-parallel, then  $K_i=\pi^{-1}_i(K_{i+1})$.
If $K$ is $E_{i+1}/E_i$-orthogonal, then the restriction $\pi_i|K_i:K_i\to K_{i+1}$ is bijective.

Now Proposition~\ref{p2.2} implies that the set $K=K_0$ is kaleidoscopic in $X=X_0$.
\end{proof}

Proposition~\ref{p2.3} suggests the following notion that will be cenral in our subsequent discussion.

\begin{definition}\label{d2.4} A (kaleidoscopic) subset $K$ in a $G$-space $X$ is called {\em splittable} if there is an increasing sequence of  $G$-invariant equivalence relations $$\Delta_X=E_0\subset E_{1}\subset\dots\subset E_m=X\times X$$ such that for every $i<m$ the set $K$ is either $E_{i+1}/E_i$-parallel or $E_{i+1}/E_i$-orthogonal.
\end{definition}

Proposition~\ref{p2.3} implies that each splittable subset in a transitive $G$-space is kaleidoscopic. What about the inverse implication?

\begin{problem} For which $G$-spaces $X$ every kaleidoscopic configuration  $K\subset X$ is splittable?
\end{problem}

\section{Kaleidoscopical configurations in matric spaces}\label{s:prim}
Here we consider each metric space $(X,d)$ as a $G$-space endowed
with the natural action of its isometry group $G=\Iso(X)$. If this
action is transitive, then the metric space $X$ is called {\em
isometrically homogeneous}.

Let us recall that a metric space $(X,d)$ is {\em ultrametric} if the metric $d$ satisfies the strong triangle inequality
$$d(x,z)\le\max\{d(x,y),d(y,z)\}$$
for all $x,y,z\in X$. It follows that for every $\e\ge 0$ the relation
$$E_\e=\{(x,y)\in X^2:d(x,y)\le\e\}\subset X\times X$$ is an invariant equivalence relation on $X$.

\begin{theorem}\label{t2.7}
Let $(X,d)$ be an isometrically homogeneous ultrametric space with the finite distance scale $d(X\times X)=\{\varepsilon_0,\varepsilon_1,\ldots,\varepsilon_n\}$ where $0=\varepsilon_0<\varepsilon_1<\ldots<\varepsilon_n$. Then every kaleidoscopical configuration $K$ in $X$ is $(E_{\varepsilon_0},E_{\varepsilon_1},\ldots,E_{\varepsilon_n})$-splittable.
\end{theorem}

\begin{proof} Assume conversely that $K$ is not $(E_{\varepsilon_0},E_{\varepsilon_1},\ldots,E_{\varepsilon_n})$-splittable.
Then for some $k<n$ the set $K$ is neither $E_{\e_{k+1}}/E_{\e_k}$-parallel nor $E_{\e_{k+1}}/E_{\e_k}$-orthogonal. We can assume that $k$ is the smallest number with that property. By $[x]_{\e_i}$ we shall denote the closed $\e_i$-ball $[x]_{E_{\e_i}}$ centered at a point $x\in X$.

Since $K$ is not $E_{\e_{k+1}}/E_{\e_k}$-orthogonal, there are two points $u,v\in K$ such that $\e_k<d(u,v)=\e_{k+1}$.
Since $K$ is not $E_{\e_{k+1}}/E_{\e_k}$-parallel, there are points $w\in K$ and $z\in X$ such that $\e_{k}<\inf_{x\in K}d(z,x)=d(z,w)=\e_{k+1}$.

Since $X$ is isometrically homogeneous, we can find an isometry
 $\varphi:X\to X$ such that $\varphi(w)=z$. Then $\varphi([w]_{\e_{k}})=[z]_{\e_{k}}$ and we can define an isometry $\phi:X\to X$ letting
$$\phi(x)=\begin{cases}
\varphi(x)&\mbox{ if $x\in[w]_{\e_{k}}$},\\
\varphi^{-1}(x)& \mbox{ if $x\in[z]_{\e_{k}}$},\\
x&\mbox{otherwise}.
\end{cases}$$
The isometry $\phi$ swaps the balls $[w]_{\e_k}$ and $[z]_{\e_k}$ but does not move points outside the union $[w]_{\e_k}\cup [z]_{\e_k}$. Since $K$ is $\chi$-kaleidoscopic, the restrictions $\chi|\phi(K)$ and $\chi|K$ are bijections onto $C$. Consequently, $\chi(w)=\chi(z')$ for some point  $z'\in[z]_{\e_k}$.
Taking into account that $d(w,z')=d(w,z)=\e_{k+1}=d(u,v)$ and $X$ is an isometrically homogeneous ultrametric space, we can construct an isometry $\psi:X\to X$ such that $\psi(u)=w$ and $\psi(v)=z'$. For this isometry, $w,z'\in\psi(K)$ and hence $\chi|\psi(K)$ is not injective, contradicting the choice of the coloring $\chi$.
\end{proof}

\begin{problem}
Let $\{0,1\}^\omega$ be the Cantor space endowed with the standard ultrametric generating the product topology. Describe all kaleidoscopical configurations in $\{0,1\}^\omega$.
\end{problem}

\begin{remark} All closed kaleidoscopical configurations in $\{0,1\}^\omega$
can be characterized with usage of Theorem \ref{t2.7}. Among them
there are plenty of non-splittable configurations.
\end{remark}

A $G$-space $X$ is called {\em primitive} if each $G$-invariant equivalence relation on $X$ is equal to $\Delta_X$ or to $X\times X$. It follows that  each splitting configuration $K$ in a primitive $G$-space $X$ is trivial, i.e. either $K=X$ or $K$ is singleton. It is natural to ask if every kaleidoscopical configuration in a primitive $G$-space trivial.

The answer to this question is affirmative if $X$ is {\em 2-transitive} in the sense that for any pairs $(x,y),(x',y')\in X^2\setminus \Delta_X$ there is $g\in X$ such that $(x',y')=(gx,gy)$.

An example of a primitive $G$-space, which is not 2-transitive is the  Euclidean space $\mathbb{R}^n$ of dimension $n\ge 2$ endowed with the action  of its isometry group $\mathrm{Iso}(\mathbb{R}^n)$. It turns out that $\IR^n$ contain $2^{\mathfrak c}$ many unsplittable kaleidoscopic configurations of cardinality $\mathfrak c$.

To construct a kaleidoscopic subset in $\IR^n$ use
proposition\ref{p1.10} and the following auxiliary definition.

Let $(X,d)$ be a metric space. By $S(x,r)=\{y\in X:d(x,y)=r\}$ we shall denote the sphere of radius $r$ centered as a point $x\in X$.

\begin{definition} A subset $K$ of a metric space $(X,d)$ is called {\em rigid} if for any pairwise distinct points $x,y,z\in K$ and numbers $r_x,r_y,r_z\in d(K\times K)$ the spheres $S(x,r_x)$, $S(y,r_y)$, $S(z,r_z)$ have no common point in $X\setminus K$.
\end{definition}

\begin{theorem}\label{t:indep} Let $X$ be metric space and $G\subset\Iso(X)$ be a group of isometries of $X$. Each infinite rigit subset $K\subset X$ of cardinality $|K|\ge|G|$ is kaleidoscopical.
\end{theorem}

\begin{proof}  The kaleidoscopicity of the set $K$ will follow from proposition~\ref{p1.10} as soon as we check that the  hypergraph $(V,\F)=(X,\{gK:g\in G\})$ satisfies the conditions (1)--(2) for the cardinal $\kappa=|K|$. Since $|G|\le \kappa=|K|=|gK|$ for all $g\in G$, the condition (1) is satisfied.

To show that (2) holds, take any subset $A\subset G$ of cardinality
$|A|<\kappa$ and any subset $B\in X\setminus AK$ of cardinality
$|B|<\kappa$. We need to show that $|\St(B,\F)\cap AK|<\kappa$. This
will follow from $\max\{|A|,|B|\}<\kappa$ as soon as we check that
$|\St(b,\F)\cap aK|\le 2$ for every $b\in B$ and $a\in A$. Assuming
conversely that $\St(b,\F)\cap aK$ contains three pairwise distinct
points $x,y,z$ we shall obtain a contradiction with the metric
independence of $K$ as $d(b,x),d(b,y),d(b,z)\in d(K\times K)$ and
$b$ is the common  point of the spheres $S(x,d(b,x))$,
$S(y,d(b,y))$, $S(z,d(b,z))$.
\end{proof}

In light of Theorem~\ref{t:indep} it is important to construct a
rigit subsets in metric spaces.
%

\begin{lemma}\label{l3.6} Any algebraic independent subset $L$ of affine line in the Euclidean space $\IR^n$ of dimension $n\ge 1$ is rigit.
\end{lemma}

\begin{proof} Identify algebraic independent $L$ with a subset of $\IR$
and let $Y$ be any subset of $\L$ with cardinality less then
$\mathfrak{c}$. It's enough to show that there are no $a\neq b\neq
c\in Y$, $r_a,r_b,r_c\in d(Y\times Y)\setminus\{0\}$ and
$x\in\IR^n\setminus(Y\cup\{p\})$ with $d(x,a)=r_a$, $d(x,b)=r_b$,
$d(x,c)=r_c$. It follows from the theorem of cosines applying to
$\cos(\angle abx)=-\cos(\angle cbx)$ and the observation that there
are no such values that
$((a-b)(b-c)+r_b^2)(a-c)-(a-b)r_c^2-(b-c)r_a^2=0\}$. Now the proof
of the last statement. Let $r_a=s_1-s_2$,$r_b=z_1-z_2$,
$r_c=t_1-t_2$ where $s_1,s_2,t_1,t_2,z_1,z_2\in A$. It is a polinom
with variable $a$ taking $r_a,r_b,r_c$ as linear functions of $a$ if
some of the numbers $s_1,s_2,t_1,t_2,z_1,z_2$ are equal to $a$ or
constants. If $z_1$ or $z_2$ is equal to $a$ then $t_1$ or $t_2$ is.
Then the coefficient of $a^2$ in the equation is $b-c-2z_2+b+2t_2$
or $b-c-2z_2+b+2t_2+b-c$ and in any case is nonzero that is
impossible. The same if one of $\{z_1,z_2\}$ is equal to $c$. If
$z_1$ or $z_2$ equals to $b$ and none of them equals to $a$ or $c$
then . If $z_1\neq a,b,c$ and $z_2\neq a,b,c$ then as a polinom with
variables $z_1,z_2$ the coefficient of $z_1z_2$ is $0$ only if
$|z_1-z_2|=|s_1-s_2|=|t_1-t_2|$ but then $(a-b)(b-c)(a-c)=0$.
\end{proof}

Now we are able to prove the promised:

\begin{theorem} Let $L$ is continual algebraic independent subset of a line in $\IR^n$. Each subset $K\subset L$ of cardinality $\mathfrak c$ is kaleidoscopic in $\IR^n$. Consequently, the Euclidean space $\IR^n$ contains $2^{\mathfrak c}$ many kaleidoscopic subsets.
\end{theorem}
\begin{proof} It follows from Theorem~\ref{t:indep} and Lemma~\ref{l3.6}
\end{proof}

\begin{problem} The Euclidean space $\IR^n$ of dimension $n\ge 2$, does it contain a non-trivial finite or countable kaleidoscopical subset $K\subset \IR^n$? {\rm If such a set $K$ exists, then its cardinality $|K|$ is not less that the chromatic number of $\IR^n$.}
\end{problem}

Let us recall that the {\em chromatic number} $\chi(X)$ of a metric space $X$ is equal to the smallest number $\kappa$ of colors for which there is a coloring of $X$ without monochrome points on the distance $1$. It is known that $4\le \chi(\IR^2)\le 7$ but the exact value of $\chi(\IR^2)$ is not known. There is a conjecture that $\chi(\IR^n)=2^{n+1}-1$, see \cite[\S47]{Soifer}.

\begin{problem} Is every finite kaleidoscopical configuration in a
(finite) primitive $G$-space trivial?
\end{problem}
Some examples of infinite $G$-spaces with only trivial finite
kaleidoscopical configurations can be found in \cite[chapter 8]{BP}

A space $\IR^n$ can also be considered as a $G$-space with respect
to the group $G=Aff(\IR^n)=\{\lambda
x+a:\lambda\in\IR\setminus\{0\},a\in\IR^n\}$ of all affine
transformations. The only kaleidoscopical configurations $K$ of
cardinality $|K|<\mathfrak{c}$ in this space are singletons as any
line that contains more then one point of kaleidoscopical
configuration has no distinct points of the same color. On the other
hand, every affine subspace of $\IR^n$ is kaleidoscopical.

\begin{question} Is there any non-splitting kaleidoscopical
configuration in $\IR^n$ with action of $Aff(\IR^n)$
\end{question}
Restricting ourself with only translations of $\IR^n$, we get a
kaleidoscopical configuration of any size $\kappa$,
$1\le\kappa\le\mathfrak{c}$. It follows from well-known
decomposition of $\IR^n$inthe direct sum of rationals and the
observation that $\IZ$ has a kaleidoscopical configuration of any
finite size.
\section{Haj\'os properties in groups and $G$-spaces}\label{s:Hajos}

In this section we reveal the relation of splittability of kaleidoscopic configurations in finite Abelian groups to the Haj\'os property introduced in \cite{Hajos} and studied in \cite{Sands62}, \cite{Sz}, \cite{SS}.

We recall that an Abelian group $G$ has the {\em Haj\'os property} if for each factorization $G=AB$ either $A$ or $B$ is periodic. A subset $A$ of a group $G$ is called {\em periodic} if $A=gA$ for some non-zero element $g\in G$.  Finite Abelian groups with Haj\'os property were classified in \cite{Sands62}:

\begin{theorem}[Haj\'os-Sands]\label{HS} A finite Abelian group $G$ has the Haj\'os property if and only if $G$ is isomorphic to a subgroup of a group that has one of the following types:

$(p^n,q)$, $(p^2,q^2)$, $(p^2,q,r)$, $(p,q,r,s)$, $(p,p)$, $(p,3,3)$, $(3^2,3)$,

$(p^3,2,2)$, $(p^2,2,2,2)$, $(p,2^2,2)$, $(p,2,2,2,2)$, $(p,q,2,2)$, $(2^n,2)$, $(2^2,2^2)$,

\noindent where $p<q<r<s$ are distinct primes and $n\in\IN$.
\end{theorem}

A group $G$ is of type $(n_1,\dots,n_k)$ if $G$ is isomorphic to the direct sum of cyclic groups $C_{n_1}{\oplus}\cdots{\oplus}\; C_{n_k}$.

Now let us define two weakenings of the Haj\'os property.

\begin{definition}\label{d3.2} An Abelian group $G$ is defined to have
\begin{itemize}
\item the {\em semi-Haj\'os property} if each complemented subset $A\subsetneq G$ either is periodic or has a periodic complementer factor in $G$;
\item the {\em demi-Haj\'os property} if for each factorization $G=A B$ one of the factors $A,B$ either is periodic or has a periodic complementer factor.
\end{itemize}
\end{definition}

It is clear that for each Abelian group $G$
$$\mbox{Haj\'os } \Ra \mbox{ semi-Haj\'os }\Ra\mbox{ demi-Haj\'os}.$$

\begin{problem} Is the semi-Haj\'os property of finite Abelian groups equivalent to the demi-Haj\'os property?
\end{problem}

The demi-Haj\'os property was (implicitly) defined in  \cite{Sands74} and follows from the quasi-periodicity of any factorization of the group.
In contrast to the Haj\'os property, at the moment we have no classification of finite Abelian groups possessing the demi-Haj\'os property. It is even not known if each finite cyclic group has the demi-Haj\'os property, see Problem 5.4 in \cite{SS}. The best known positive result on the semi-Haj\'os property is the following version of Theorem 5.13 \cite{SS}:

\begin{theorem}[Bruijn-Szab\'o-Sands]\label{semiH} Each finite Abelian group $G$ of square-free order $|G|$ has the semi-Haj\'os property.
\end{theorem}

We say that a number $n$ is {\em square-free} if $n$ is not divisible by the square $p^2$ of any prime number $p$.

Surprisingly, the following problem of Fuchs and Sands
\cite[p.364]{Fuchs}, \cite{Sands74}, \cite[p.120]{SS} posed in 60-ies still is open:

\begin{problem} Has each finite Abelian group the demi-Haj\'os property?
\end{problem}

The ``semi'' version of this problem also is open:

\begin{problem}\label{pr3.4} Has each finite Abelian group the semi-Haj\'os property?
\end{problem}

The semi-Haj\'os property is tightly connected with the splittability of kaleidoscopical configurations. In order to state the precise result, let us generalize the definition of the semi-Haj\'os property to $G$-spaces.

\begin{definition} A $G$-space $X$ has the {\em semi-Haj\'os property} if  for each kaleidoscopic subset $K\subsetneq X$ there is a $G$-invariant equivalence relation $E\ne\Delta_X$ on $X$ such that $K$ is $E$-parallel or $E$-orthogonal and the set $K/E$ is kaleidoscopic in the $G$-space $X/E$.
\end{definition}

For finite Abelian groups this definition of the semi-Haj\'os property agrees with that given in Definition~\ref{d3.2}.

\begin{proposition} A finite Abelian group $G$ has the semi-Haj\'os property if and only if it has that property as a $G$-space.
\end{proposition}

\begin{proof} Assume that the group $G$ has the semi-Haj\'os property.
To show that the $G$-space $G$ has the semi-Haj\'os property, take any kaleidoscopic subset $A\subset G$. By Corollary~\ref{c1.2}, $A$ is complementable and hence has a complementer factor $B$. Since $G$ has the semi-Haj\'os property, either $A$ is periodic or else $A$ has a periodic complementer factor. In the latter case we can assume that the complementer factor $B$ is periodic. Consequently there is a non-trivial cyclic subgroup $H\subset G$ such that either $A+H=A$ or $B+H=B$.
Consider the quotient group $G/H$ and the quotient homomorphism $q:G\to G/H$. By Lemma 2.6 of~\cite{SS}, the images $A/H=q(A)$ and $B/H=q(B)$ form a factorization $G/H=A/H\cdot B/H$ of the quotient group $G/H$. Consequently, the set $A/H$ is complemented in $G/H$ and by Corollary~\ref{c1.2}, it is kaleidoscopic in $G/H$.

The subgroup $H$ induces a $G$-invariant equivalence relation $E=\{(x,y)\in G:x-y\in H\}$ whose quotient space $G/E$ coincides with the quotient group $G/H$. We claim that the set $A$ is either $E$-parallel or $E$-orthogonal. By the choice of the group $H$, we get $A=A+H$ or $B=B+H$. In the first case the set $A$ is $E$-parallel. In the second case $A$ is $E$-orthogonal as $(A-A)\cap H\subset (A-A)\cap (B-B)=\{0\}$.
\smallskip

Now assuming that the $G$-space $G$ has the semi-Haj\'os property, we shall prove that the group $G$ has the semi-Haj\'os property. Given any complemented subset $A\subset G$ we need to show that either $A$ is periodic or else $A$ has a periodic complementer factor. By
Corollary~\ref{c1.2}, the set $A$ is kaleidoscopic in the $G$-space $G$. The semi-Haj\'os property of the $G$-space $G$ guarantees the existence of an invariant equivalence relation $E\ne\Delta_G$ on $G$ such that $A$ is $E$-parallel or $E$-orthogonal and $A/E$ is kaleidoscopic in $G/E$. It follows that the equivalence class $H=[0]_E$ of zero is a subgroup of the group $G$. Taking into account that $E$ is $G$-invariant, we conclude that $(x,y)\in E$ iff $x-y\in[0]_E$. So, $G/E$ coincides with the quotient group $G/H$. The set $A/H$, being kaleidoscopical, is complemented in $G/H$ according to Corollary~\ref{c1.2}. Consequently, there is a subset $B_H\subset G/H$ such that $G/H=A/H\cdot B_H$. Let $q:G\to G/H$ be the quotient map and $s:G/H\to G$ be any section of $q$.

Now consider two cases. If $A$ is $E$-parallel, then $A=A+H$ is periodic and complemented as $B=s(B_H)$ is a complementer factor to $A$ in $G$. If $A$ is $E$-orthogonal, then the complete preimage $B=q^{-1}(B_H)$ is a periodic complementer factor to $A$ in $G$.
\end{proof}

Now we reveal the relation between the semi-Haj\'os property and the splittability of kaleidoscopic sets.

\begin{proposition}\label{p3.7} If each kaleidoscopic subset of a transitive $G$-space $X$ is splittable, then $X$ has the semi-Haj\'os property.
\end{proposition}

\begin{proof} To show that $X$ has the semi-Haj\'os property, fix any kaleidoscopic subset $K\subset X$. By our assumption, $K$ is $(E_0,\dots,E_m)$-splittable by some increasing chain of invariant equivalence relations $\Delta_X=E_0\subset \cdots\subset E_m=X\times X$.  For every $i\le m$ consider the quotient $G$-space $X_i=X/E_i$ and let $q_i:X\to X_i$ be the quotient projection. Also let $K_i=q_i(K)\subset X_i$. By Proposition~\ref{p2.2}, $K_i$ is kaleidoscopic in the $G$-space $X_i$.
In particular, $K_1$ is kaleidoscopic in $X_1=X/E_1$. By Definition~\ref{d2.4}, $K=K_0$ is either $E_1$-parallel or $E_1$-orthogonal. This means that $X$ has the semi-Haj\'os property.
\end{proof}

Theorems~\ref{t2.7} and Proposition~\ref{p3.7} imply

\begin{corollary} Each isometrically homogeneous ultrametric space with finite distance scale has the semi-Haj\'os property.
\end{corollary}

A $G$-space $Y$ is defined to be a {\em quotient} of a $G$-space $X$ if  $Y$ is the image of $X$ under a $G$-equivariant map $f:X\to Y$.

\begin{proposition}\label{p3.9} Each kaleidoscopical subset of a $G$-space $X$ is splittable provided that:
\begin{enumerate}
\item each quotient $G$-space of $X$ has the semi-Haj\'os property and
\item $X$ admits no strictly increasing infinite sequence $(E_n)_{n\in\w}$ of $G$-invariant equivalence relations.
\end{enumerate}
\end{proposition}

\begin{proof} Assume that some kaleidoscopic subset $K\subset X$ is not splittable. Let $K_0=K$, $E_0=\Delta_X$, and $X_0=X/E_0=X$. Since $X$ has the semi-Haj\'os property, there is a $G$-invariant equivalence relation
$E_1\ne \Delta_X$ on $X_0$ such that the set $K_1=K_0/E_0$ is kaleidoscopic in the $G$-space $X_1=X_0/E_1$ and $K_0$ is either $E_1$-parallel or $E_1$-orthogonal.

By our assumption, $K$ is not splittable, so $X_1$ is not a singleton. The $G$-space $X_1=X/E_1$, being a quotient of $X$, has the semi-Haj\'os property. Consequently, for the kaleidoscopic set $K_1\subset X_1$ there is a $G$-invariant equivalence relation $\tilde E_2\ne\Delta_{X_1}$ on $X_1$ such that the set $K_1$ is $\tilde E_2$-parallel or $\tilde E_2$-orthogonal and the quotient set $K_2=K_1/\tilde E_1$ is kaleidoscopical in the $G$-space $X_2=X_1/\tilde E_1$. Let $q^1_2:X_1\to X_2$ be the quotient projection. The composition $q^1_2\circ q_1:X\to X_2$ determines the $G$-invariant equivalence relation $E_2=\{(x,x')\in X^2:q^1_2\circ q_1(x)=q^1_2\circ q_1(x')\}$ on $X$ such that $X/E_2=X_2$ and $K_2=K/E_2$ and $K_1$ is either $E_2/E_1$-parallel or $E_2/E_1$-orthogonal.

Continuing by induction, we shall produce an infinite increasing sequence $(E_n)_{n\in\w}$ of $G$-invariant equivalence relations on $X$ such that for every $n\in\IN$ the set $K_n=K/E_n$ is kaleidoscopic in the $G$-space $X/E_n$
and $K$ is either $E_{n}/E_{n-1}$-parallel or $E_n/E_{n-1}$-orthogonal.
But the existence of an infinite strictly increasing sequence of $G$-invariant equivalence relations on $X$ contradicts our assumption.
\end{proof}

Since each quotient group of a finite Abelian group $G$
is isomorphic to a subgroup of $G$, Proposition~\ref{p3.9} implies:

\begin{corollary}\label{c3.10} If each subgroup of a finite Abelian group $G$ has the semi-Haj\'os property, then each kaleidoscopic subset $K\subset G$ is  splittable.
\end{corollary}

\begin{question} Assume that a finite Abelian group $G$ has  the semi-Haj\'os property. Has each subgroup of $G$ that property?
\end{question}

The classification of finite Abelian groups with Haj\'os property
given in Theorem~\ref{HS} implies that this property is inherited by subgroups. Because of that, Corollary~\ref{c3.10} implies:

\begin{corollary} For a finite Abelian group $G$ with the Haj\'os property, each kaleidoscopical subset $K\subset G$ is splittable.
\end{corollary}

Also Proposition~\ref{p3.9} and Theorem~\ref{semiH} imply:

\begin{corollary} For a finite Abelian group $G$ of square-free order $|G|$ each kaleidoscopical subset $K\subset G$ is splittable.
\end{corollary}

\begin{remark} It follows from Proposition~\ref{p3.7} and Corollary~\ref{c3.10} that Problems~\ref{pr2.6} and \ref{pr3.4} are equivalent (and both are open and apparently difficult).
\end{remark}

According to an old result of Haj\'os \cite{Hajos}, if in a factorization  $\IZ=A+B$ of the infinite cyclic group $\IZ$ the factor $A$ is finite, then the factor $B$ is periodic. We do not know if the same is true for the groups $\IZ^n$ with $n\ge 2$.

\begin{problem} Assume that $\IZ^n=A+B$ is a factorization with finite factor $A$. Is the factor $B$ periodic? Has $A$ a periodic complementer factor?
\end{problem}


\end{document}